\documentclass[11pt]{article}
\usepackage{amsthm} 
\usepackage{amsmath} 
\usepackage{amssymb}
\usepackage{amsfonts} 
\usepackage{amscd}
\usepackage[mathscr]{eucal} 
\usepackage{verbatim}
\usepackage{graphics}

\newtheorem{thm}{Theorem}[section]  
\newtheorem*{mlem*}{Main Lemma}
\newtheorem*{thm*}{Theorem}  
\newtheorem{prop}[thm]{Proposition}
\newtheorem{lem}[thm]{Lemma}
\newtheorem*{cor*}{Corollary}
\theoremstyle{definition}
\newtheorem{defi}[thm]{Definition}
\newtheorem{rem}[thm]{Remark}
\newtheorem{Q}[thm]{Question}
\newtheorem*{conj}{Conjecture}

\newcommand{\Z}{{\mathbb Z}}

\begin{document} 

\title{Large sum-free sets in abelian groups }
\author{R.Balasubramanian, \\ The Institute of
Mathematical Sciences\\ CIT Campus, Taramani\\ Chennai-600113,
India.\\ E-mail: balu@imsc.res.in\\
\\
Gyan Prakash\\
Harish Chandra Research Institute\\
Chatnag Road, Jhusi,\\
Allahabad - 211 019, India\\
E mail: gyan@mri.ernet.in  }
\date{}
\maketitle
\begin{abstract}
Let $A$ be a subset of a finite abelian group $G$. We say that $A$ is sum-free
if the equation $x + y = z$, has no solution $(x,y,z)$ with $x, y, z$ belonging
to the set $A$. In this paper we shall characterise the largest possible
sum-free subsets of $G$ in case the order of $G$ is only divisible by primes 
which are congruent to $1$ modulo $3$.
\end{abstract}
\section{Introduction and statements of result}\label{sintro}
 Throughout this paper
$G$ will denote a finite abelian group of order $n$. If $A$ is a
subset of $G$ then we say that $A$ is sum-free if the equation 
 $ x + y = z$ has no solution with $(x, y, z) \in A \times A\times A$. We say that $A$ is maximal sum-free
if it is not a proper subset of any sum-free set.
The present article is motivated by the following question:
\begin{Q}
 Find a ``{\it structure}'' of all ``{\it large}'' maximal sum-free 
subsets of $G$.
\end{Q}
In this regard we prove  Theorem~\ref{exam}, 
Theorem~\ref{BGCH} and Theorem~\ref{maxch}.  
The results are based on Theorem~\ref{GRM} which is a
 recent result of Ben~ 
Green and Imre Ruzsa~\cite{GR}. Our results give complete characterisation
of largest sum-free subsets of $G$ in case all the divisors of $n$ are 
congruent to $1$ modulo $3$. For all other groups a structure of all largest
sum-free subsets was known before.\\
\\
Before stating our results and previously known results related to the above 
question we shall explain what do we mean by ``{\it large}'' and what sort
of '{\it structure}'' does one may expect?\\
\\
To understand the meaning of {\it large}
in the above question, we need to understand the following 
question.
\begin{Q} 
 How big is the largest
sum-free subset of $G$ ?
\end{Q}
\noindent
\newpage
\begin{defi}
\begin{enumerate}
\item Given any finite abelian group $K$ and
a set $A \subset K$, we define the density
of the set $A$ to be $\frac{|A|}{|K|}$ and denote it by $\alpha(A)$.
\item We use $\mu(G)$ to denote the density of a largest sum-free subset 
of $G$.
\end{enumerate}
\end{defi}
We say a sum-free set $A \subset G$ is large if the density of the set $A$ is 
close to $\mu(G)$; that is $\mu(G) - \alpha(A)$ is  
``{\it small}''.\\
The value of $\mu(G)$ is now known for any finite abelian group $G$ \cite{GR}.
\begin{thm}\textup{(\cite{GR}, Theorem 2)} Let $G$ be a finite abelian group and $m$ be
its exponent. \label{mug}
Then the value of $\mu(G)$ is given by the following formula.
\begin{equation*} 
\mu(G) = max_{d | m}  \frac{ \left[\frac{d -2}{3} \right] + 1}{d}
\end{equation*}
\end{thm} 
The following facts are straightforward to check.
\begin{enumerate}
\item For any positive integer $d$, we have the natural projection 
$p_d : \Z \to \Z/d\Z$. Let the set $B_d \subset \Z/d\Z$ be the image of 
integers
in the interval $(\frac{d}{3}, \frac{2d}{3}]$ under the map $p_d$. Then it is 
straightforward to check that $B_d$ is sum-free and density of the set $B_d$
is given by
\begin{equation*}
\alpha(B_d) = \frac{ \left[\frac{d -2}{3} \right] + 1}{d}. 
\end{equation*}
\item  For any positive integer $d$, there is a surjective 
homomorphism $f: G \to \Z/d\Z$ if and only if $d$ divides the exponent of the
group $G$.
\item For any positive integer $d$ and a surjective homomorphism  
$f: G \to \Z/d\Z$, the set $A = f^{-1}(B_d)$ is  a sum-free subset of $G$ and
$\alpha(A) = \alpha(B_d)$.
\end{enumerate}
Therefore the following result follows.
\begin{thm}(\cite{GR}) \label{struct}
 Given any finite abelian group $G$, there exists a sum-free
set $A \subset G$,  a surjective homomorphism $f: G \to \Z/d\Z$ (where $d$ is 
a positive integer), a sum-free set $B \subset \Z/d\Z$; such that the 
following hold.
\begin{enumerate}
\item $A = f^{-1}(B)$.
\item The density of $A$ is $\mu(G)$.
\end{enumerate} 
\end{thm}
Now regarding the structure of {\it large} sum-free subsets of $G$ we may ask
the following question.
\begin{Q}
 Let $A$ be a ``large'' sum-free subset of $G$.
Then given any such  set $A$, does there always
 exist a surjective homomorphism 
$f: G \to \Z/d\Z$ (where $d$ is some positive integer)  
and a sum-free set $B \subset
\Z/d\Z$, such that the set $A$ is a subset of the set $f^{-1}(B)$?
\end{Q}
Before discussing this question we 
describe the value of  $\mu(G)$ more explicitly by dividing the finite
abelian groups in the following three classes.
\begin{defi}
 Suppose that $G$ is a finite abelian group of order
$n$. 
\begin{enumerate}
\item If  $n$ is divisible by any prime $p \equiv 2 {\rm (mod \ 3)}$
then we say that $G$  is type $I$. We say that G is type $I(p)$ if
it is type $I$ and if $p$ is the {\it least} prime factor of $n$ of the form
$3l + 2$.  In this case the value of $\mu(G)$ is equal to 
$\frac{1}{3}~+~\frac{1}{3p}$.
\item If $n$ is not divisible by any
prime $p \equiv 2 {\rm (mod \ 3)}$,  but $3|n$, then we say that $G$
is type $II$. In this case the value of $\mu(G)$ is equal to $\frac{1}{3}$.
 \item The group $G$
is said to be of type $III$ if all the divisor of $n$ (order of $G$) are 
congruent to 
$1$ modulo $3$. Let $m$ be the exponent of $G$. In this case  
the value of $\mu(G)$ is equal to 
$\frac{1}{3}~-~\frac{1}{3m}$. We also  note the fact that if $G$ is a 
group of type $III$ then any subgroup as well as quotient of $G$ is also a
type $III$ group.
\end{enumerate}
\end{defi}
\begin{rem}
We note the fact that if $G$ is a type $III$ group and $d$ is any divisor
of $m$, then $d$ is odd and congruent to $1$ modulo $3$. Therefore  
$d$ is congruent to $1$ modulo $6$ and $\frac{d -1}{6}$ is a non~negative
integer.
\end{rem}
Theorem~\ref{mug} was proven for type $I$ and
type $II$ groups by 
Diananda and Yap~ \cite{DY}. For some special cases of type $ III $ groups it 
was proven by various authors
( see \cite{Yap1, Yap, RS} ). For an arbitrary abelian groups of 
type $III$ the proof of Theorem~\ref{mug} is due to 
Ben Green and Ruzsa~\cite{GR}.\\
\\
\noindent
Hamidoune and
Plagne~\cite{Plagne1} answered the question $3$ affirmatively when 
$|A| \geq \frac{|G|}{3}$, in the case $|G|$ is odd. In case $|G|$ is even
they answered the question $3$ affirmatively if $|A| \geq \frac{|G| + 1}{3}$.
In case  $G
= (\Z/2\Z)^r$ with $r \geq 4$ and $|A| \geq 5. 2^{r - 4}$ then Davydov
and  Tombak~\cite{Davy} showed that  answer of this question is
affirmative.  Recently Lev~\cite{largelev} answered  this question
affirmatively in the case when $G = (\Z/3\Z)^r$ (with an integer
$r \geq 3$) and $|A| >\frac{ 5}{27} 3^{r}$. Lev~\cite{Levp} has also
 characterised
the sum-free subsets $A$ of $ \Z/p\Z$ when $p$ is prime and  $|A| >
0.33p$.\\
\\
\noindent
Notice that none of the above mentioned results tells us anything related
to the question~3, in case
$G$ is a finite abelian group of type $III$ and  $G$ is not cyclic.
In case $G$ is cyclic the answer of question~$3$ is obviously affirmative.
In case $G$ is not cyclic and of type $III$ Theorem~\ref{exam} shows that
the answer of  question~$3$ is negative. \\ \\
\noindent
For the rest of this paper unless specified differently, $G$ shall denote a 
finite
abelian group of type $III$ and of order $n$.
The symbol $m$ shall denote the exponent of $G$
and $k = \frac{m -1}{6}$. 
\begin{rem}
If $G$ is an abelian group of type $III$ and $m$ is exponent of $G$, 
there exist $S \subset G$ and $C \subset G$ such that $S$ and $C$ are 
 subgroups of $G$, $C$ is isomorphic to $\Z/m\Z$ 
and $G = S \oplus C$. In  case $G$ is not cyclic, $S$ will be a nontrivial
subgroup of $G$.
\end{rem}
Let $p: \Z \to \Z/m\Z$ be the natural projection 
from the group of integers to $\Z/m\Z$. 
\begin{thm} \label{exam} Let $G$ be a finite abelian group of type
$III$. Let $m$ denote the exponent of $G$ and  $k = \frac{m -1}{6}$.
Suppose $G$ is not a cyclic group and $S, C$ are non trivial 
subgroups of $G$ such that
 $G =
S \oplus C$,  and $C$ is isomorphic to
$\Z/m\Z$.
Let $g: C \to \Z/mZ$ be an isomorphism.  Let
$J$  be any proper subgroup of $S$ and $b$ is any element
belonging to $S$.  Then consider the following  two examples:
\begin{enumerate}
\item The set $A = ((J + b) \oplus g^{-1}p(\{2k\})) \cup (S \oplus g^{-1}
p(\{2k + j: 1 \le j \le 2k -1  \} ) \cup ((J +2 b)^c \oplus g^{-1}p(\{4k\}))$. (Here for any set $D \subset S$ the symbol $D^c$ denotes the
set $S \setminus D$.)
\item  The set $A = ((J + b) \oplus g^{-1}p(\{2k\})) \cup ((J- 2b)^c \oplus
g^{-1}p(\{2k + 1\}) \cup (S \oplus g^{-1}p(\{2k + j: 2 \le j \le 2k -1\}) ) 
\cup
((J + 2b)^c \oplus~g^{-1}p(\{4k\})) \cup 
((J -b)  \oplus g^{-1}p(\{4k + 1 \}))$.
\end{enumerate} 
Let $A$ be any of the set as above. Then the following holds.
\begin{enumerate}
\item The set $A$ is a sum-free subset of $G$ and $\alpha(A) =  \mu(G)$.
\item For any positive integer $d$, there does not exist any surjective 
homomorphism \\
$f:~G~\to~
\Z/dZ$, a sum-free set $B \subset Z/d\Z$, such that the set
 $A$ is a subset of the set $f^{-1}(B)$.
\end{enumerate}
\end{thm}
\noindent 
We got to know the above examples from certain remarks made in~\cite{GR} about 
the group $\left(\Z/7\Z \right)^r$. \\ \\
\noindent We prove that if $G$ is a type $III$  group and $A$ 
is a sum-free subset $G$ of largest possible cardinality then either
$A$ is an inverse image of a sum-free subset of a cyclic quotient of $G$
or $A$ is one of the set as given in Theorem~\ref{exam}. 
\begin{defi}
 Given a sum-free set $A \subset G$ ,
a surjective homomorphism $h: G \to \Z/dZ$ the following definition
and notations are useful.
\begin{enumerate}
\item For any $i \in \Z/d\Z$ the symbol $A(h, i)$ denote the set $A \cap
h^{-1}\{i\}$.
\item For any
$i \in \Z/d\Z$
we define $\alpha(h, i)
= \frac{d}{n}|A(h, i)|$.
\item Let $l = \frac{d -1}{6}$ and $p_d: \Z \to \Z/d\Z$ be the natural 
projection. The sets $H_d
,T_d, M_d, I_d \subset \Z/d\Z$ denote the  following sets.\\
\begin{eqnarray*} H_d &=& p_d\{l + 1, l + 2, \cdots, 2l -1, 2l \}\\   M_d &=&
p_d\{2l + 1, 2l + 2, \cdots, 4l -1, 4l \}\\ T_d &=& p_d\{4l + 1, 4l + 2,
\cdots, 5l -1, 5l \}\\ I_d &=& p_d\{l + 1, l + 2, \cdots, 5l -1, 5l \}\\
\end{eqnarray*}
\item The symbol $H
,T, M, I$ denote the sets $H_m, T_m, M_m, I_m$ respectively. The symbol $p$ 
denotes the map $p_m$.  
\end{enumerate}
\end{defi}
\begin{thm}\label{BGCH} 
Let $G$ be a finite abelian group of type $III$. Let $A$ be a sum-free 
subset of an abelian
group $G$. Let $n$ denotes the order of $G$ and $m$ denotes the exponent of
$G$.
Let $k = \frac{ m -1}{6}$ and
$\eta = 2^{-23}$.  Suppose that  $\alpha(A) > \mu(G) - min(\eta,
\frac{5}{42m})$. Then there exists a surjective homomorphism
\begin{equation*} f: G \to \Z/m\Z
\end{equation*} such that the following holds. 
\begin{equation*}
A \subset f^{-1}p(\{2k + j: 0 \le j \le 2k + 1  \})
\end{equation*}
Further the following also holds.
\begin{enumerate} 
\item For all 
$i \in p(\{2k + j: 2 \le j \le 2k -1 \})$ 
the inequality $ \alpha(f,i) \geq~1 - m \left(\mu(G) - 
\alpha(A)  \right)$ holds.
\item $\alpha(f, 2k) + \alpha(f, 4k) \geq 
~1 - m \left(\mu(G) - 
\alpha(A)  \right)$.
\item $\alpha(f, 4k + 1) + \alpha( f, 2k + 1) \geq 
~1 - m \left(\mu(G) - 
\alpha(A) \right)$
\item If $A$ is maximal then 
\[
f^{-1}p(\{2k + j: 2 \le j \le 2k -1  \})  \subset A
\]
\end{enumerate}
\end{thm}
\noindent
Using Theorem~\ref{BGCH} the following theorem follows easily.
\begin{thm}\label{maxch} Let $A$ be a sum-free subset of an abelian
group $G$ of type $III$. Let  the symbol $m$ denote the exponent of $G$ and
 $ k = \frac{m-1}{6}$. Let the density of the  set $A$ is equal to 
$ \mu(G)$ and 
$f: G \to \Z/m\Z$ be a surjective homomorphism as given by Theorem~\ref{BGCH}.
 Let the set $S$ denotes the kernel of $f$ and $C$ is 
a subgroup of $G$ such that  $G = S \oplus C$. Let $g: C \to \Z/m\Z$ be 
an isomorphism obtained by restricting $f$ to the subgroup  $C$.
Then there exist $J$
a subgroup of $S$ and $b \in S$ such that one of the
following  holds:
\begin{enumerate}
\item The set $A = f^{-1} p\left(\{2k + j : 1 \le j 
\le 2k \} \right)$.
\item One of the set $A$ or $-A$ is equal to the following set 
\begin{equation*} 
\left((J + b) \oplus  g^{-1}p(\{2k\})\right) \cup 
\left(f^{-1} p(\{2k + j: 1 \le j \le 2k -1\} ) \right) 
\cup ((J + 2b)^c \oplus g^{-1} \{4k\}).
\end{equation*} (Here and
in the following for any set $D \subset S$ the symbol $D^c$ denotes the
set $S \setminus D$.)
\item  The set $A$ is union of the sets
$f^{-1} 
p(\{2k + j: 2 \le j \le 2k -1 \} )$, 
$\left((J + b)~\oplus g^{-1}p(\{2k\})\right)$, 
$\left((J + 2b)^c~\oplus~g^{-1}p(\{4k\})\right)$,
  $\left((J~-~b)~\oplus~g^{-1} p(\{4k~ +~ 1 \})\right)$ \\ and
$\left((J- 2b)^c~\oplus~g^{-1}p(\{2k + 1\})\right)$.
\end{enumerate}
\end{thm}
\noindent
The proof of Theorem~\ref{BGCH} is based on the following
result of Ben Green  and Ruzsa~\cite{GR}.
\begin{thm} \label{GRM} (\cite{GR}, Proposition 7.2.)
 Let $A$ be a sum-free subset of an abelian
group $G$ of type $III$. Let $\eta = 2^{-23}$. 
Suppose that  $\alpha(A) \geq \mu(G) -\eta$, then there exists a surjective
homomorphism $\gamma:G \to \Z/q\Z$ with $q \neq 1$ such  that the
following holds.
\begin{equation*}
A \subset
\gamma^{-1}I_q.
\end{equation*}
\end{thm}
\noindent 
We require the following definitions and notations.
\begin{defi}
 Let $A$ be a sum-free subset of $G$ and
$\alpha(A)  \geq \mu(G) -\eta$, then we choose a $\gamma$ a surjective 
homomorphism from $Z$ to $\Z/q\Z$ with $q \ne 1$ and 
$A  \subset
\gamma^{-1}I_q$. Following the terminology
of~\cite{GR} we call $\gamma$ the special direction of the set $A$
and $q$ the order of special direction. We use the symbols $\alpha_i$ and
$A_i$ to denote
 the number $\alpha(\gamma, i)$ and the set $A(\gamma, i)$ respectively.
We use symbol $S$ 
to denote the 
set $ker(\gamma)$ and $S_i$ to denote $\gamma^{-1}\{i\}$. 
\end{defi}
\begin{rem}
 Notice that in case $G = \left(\Z/7\Z
\right)^r$, the  Theorem~\ref{GRM} is equivalent to
Theorem~\ref{BGCH}, as in this special case $l + 1 = 2l$.
\end{rem}
\section{An outline of the paper}
Let $H(A)$ is the largest subset of $G$ such that $H(A) + A = A$.
 The set $H(A)$ as defined is called period or stabiliser of the set $A$.
For any set $A$ as as in Theorem~\ref{exam} we prove that $H(A) = J$ where
$J$ is as in Theorem~\ref{exam}. 
Using this Theorem~\ref{exam} is easy to prove.\\
\\
\noindent 
Let $A$ be as in Theorem~\ref{BGCH} and $\gamma$ be the special direction
of $A$.
The proof of Theorem~\ref{BGCH} is divided into following three parts. 
\begin{enumerate}
\item The order of the special direction of  $A$ is $m$.
\item We define the set $L = \{i \in \Z/m\Z: \alpha_i > \frac{1}{2}\}$.
We show that for any $i, j \in L$, $\alpha_{i + j = 0}$. In particular
 the set $L$ is a sum-free subset of $\Z/m\Z$. Moreover the cardinality of 
$L$ is $2k = 2\frac{m-1}{6}$; that is $L$ is a largest sum-free subset of 
$\Z/m\Z$.
\item We describe all the largest sum-free subset of $\Z/m\Z$ which are 
subsets of $I$. Using Theorem~\ref{GRM} this characterises all the largest 
sum-free subsets of $\Z/m\Z$.
\end{enumerate}
There are two facts which we use repeatedly. One is that for any 
$i, j \in Z/q\Z$ the
sets $A(\gamma,i) + A(\gamma, j)$ and $A(\gamma, i + j)$ are disjoint.
Another is that for any divisor $d$ of $m$, the set $I_d$ is divided into
$2\frac{d -1}{6}$ disjoint pairs of the form $(i, 2i)$ with $i$ belonging to 
the set $H_d \cup M_d$.  
\section{ Stabiliser of largest sum-free subset } \label{secstab}
In this section we shall give the proof of Theorem~\ref{exam}.
Any abelian group acts on itself by translation.
Given any set $B \subset G$ we define the set $H(B)$ to be
those elements of the group $G$ such that the set $B$ is stable under the 
translation by the elements of the set $H(B)$.
In other words the set
$H(B) = \{ g \in G : g + B
= B\}$. For any set $B$ the set $H(B)$ is a subgroup.\\
\\
Let $G$ be an abelian group of type $III$ and let $A \subset G$ be as in 
Theorem~\ref{exam}.
To prove Theorem~\ref{exam} we shall prove the following 
\begin{lem}
Let $S$ and $C$ are as in Theorem~\ref{exam} and $\pi_C : G \to G/S = C$ be
the natural projection. Then the set $\pi_C(H(A)) = \{0\}$.
\end{lem}
\begin{proof}
Since $\pi_C$ is a homomorphism and $H(A)$ is a subgroup of $G$, 
 the set $\pi_C (H(A))$ is a subgroup of M. 
Since $H(A) + A = A$ by the definition of $H(A)$, it follows that
$\pi_C(H(A)) + \pi_C(A) = \pi_C(A)$. Therefore the set $\pi_C(A)$ is a 
union of cosets of $\pi_C\left(H(A)\right)$. Therefore the cardinality of the subgroup 
$\pi_C(H(A))$ divides the cardinality of the set $\pi_C(A)$.
Since the set $\pi_C(H(A))$ is a subgroup of $C$, it is also true that
 $\pi_C(H(A))$ divides
$m$. This implies that  $|\pi_C(H(A))|$ divides $gcd(|\pi_C(A)|, m)$.
Now if $A$ is a set as in Theorem~\ref{exam} $(I)$ then the cardinality of the
set $(\pi_C(A)$ is equal to $2k + 1$ and if  
$A$ is a set as in Theorem~\ref{exam} $(II)$ then the cardinality of the
set $(\pi_C(A)$ is equal to $2k + 2$. 
Since $m = 6k + 1$ it follows that in first case the number
$gcd(|\pi_C(A)|, m)$ divides 
$2$ and in the second case it divides $5$. But as $G$ is type $III$
group, any divisor of $m$  which is not equal to $1$ is
greater than or equal to $7$. Hence
$gcd(|\pi_C(A)|, m) = 1$ for any of the set $A$ as in Theorem~\ref{exam}.
This forces $|\pi_C(H(A))| = 1$ and hence the lemma follows.
\end{proof}
\begin{prop}\label{stab}
Let $A$ be any of the set as in Theorem~\ref{exam} and $S$ be a subgroup of
$G$, $J$ be a proper subgroup of $S$ 
as in Theorem~\ref{exam}. Then the stabiliser of the
set $A$ is equal to the set $J$.
\end{prop} 
\begin{proof} 
Using the previous lemma it follows that $H(A) +  \left((J + b) \oplus 
g^{-1} \{2k\} \right) = (J + b) \oplus g^{-1}\{2k\} $. This implies that 
$H(A) + J + b = J + b$. This implies that $H(A) \subset J$. But it is
straightforward to check that $J + A = A$. Therefore it follows that
 $J = H (A)$, proving the claim.
\end{proof}
\begin{proof} {\it of Theorem~\ref{exam}:}
\begin{enumerate}
\item This is straightforward to check.
\item Suppose $A$ is one of a set as in Theorem~\ref{exam}. Suppose the 
claim is
not true for this set. Then there exist a positive integer $q$, 
a surjective homomorphism
$f: G \to \Z/q\Z$, a sum-free set $B \subset \Z/q\Z$ such that the set
$A$ is a subset of $f^{-1}(B)$. Since from $(I)$, the set $A$ is a sum-free
set of largest possible cardinality, it follows that the set 
$A = f^{-1}(B)$. Therefore the kernel of $f$ is contained in the set $H(A)$.
Therefore we have the following inequality
\begin{equation}\label{g}
|H(A)| \geq \frac{n}{q}.
\end{equation} 
But from Proposition~\ref{stab} the stabiliser
of the set $A$ is $J$ which is a proper subgroup of $S$. 
Since $m$ is the exponent of $G$ it follows that $q$
is less than or equal to $m$. Therefore we have the following inequality
\begin{equation}\label{l}
|H(A)| < |S| = \frac{n}{m} \leq \frac{n}{q}
\end{equation}
This contradiction proves the claim.
\end{enumerate}
\end{proof}

\section{Order of the special direction}\label{sord}
Let $A$ be as sum-free subset of $G$ and $\alpha(A) >
min(\eta, \frac{5}{42m})$.  .
 Let $\gamma$ be the special direction of the set $A$
as given by Theorem~\ref{GRM}. In this section we shall show that the order of
$\gamma$ is equal to $m$. The proof of this result is inherent in~\cite{GR}.
We reproduce the proof here for the sake of completeness.
\begin{lem} \label{alpha1} (\cite{GR}, Lemma 7.3. (ii) )
Let $A$ be a sum-free subset of the group $G$. Let $g$ be any surjective 
homomorphism $g: G \to \Z/d\Z$, where $d$ is a positive integer.
Then for any $i \in \Z/d\Z$, the following inequality holds.
\begin{equation}
\alpha(g, i) + \alpha(g, 2i) \leq 1
\end{equation}
Here $\alpha(g, i)$ is a number as defined in section~\ref{sintro}.
\end{lem}
\begin{proof}
For any $i \in \Z/d\Z$, let the set $A(g, i)$ be as defined in 
section~\ref{sintro}. The fact that $g$ is a homomorphism implies that the set
$A(g, i) + A(g, i)$ is a subset of the set $g^{-1}\{2i\}$. The fact that the 
set
$A$ is sum-free implies that the set $A(g, i) + A(g, i)$
is disjoint from the set $A(g, 2i)$.
Therefore we have the following inequality.
\begin{equation}
|A(g, i) + A(g, i)| + |A(g, 2i)| \leq |g^{-1}\{2i\}| \label{disj}
\end{equation}  
The claim follows by observing that the set $A(g, i) + A(g, i)$ has
cardinality at~least $|A(g, i)|$.
\end{proof} 
The following lemma is straightforward to check, but is very useful.
\begin{lem}\label{2l}
Let $ d$ be a positive integer congruent to $1$ modulo $6$. 
Let $d = 6l + 1$. Let the set
$I_d, H_d, M_d, T_d$ are subsets of the group $\Z/d\Z$ as defined in 
section~\ref{sintro}.  The set $I_d$ is divided into $2l$ disjoint pairs of the
form $(i, 2i)$ where $i$ belong to the set $H_d \cup T_d$ and $2i$ belong 
to the set $M_d$.
\end{lem}
\begin{prop}\label{ord}  Let $A$ be a sum-free subset of an abelian group
 $G,$ of type $III$. Let $m$ be the exponent of $G$ and $\alpha(A)
> min(\eta, \frac{5}{42m})$. Let
$\gamma: G \to \Z/q\Z$ be  the special direction of the set  as given by
Theorem~\ref{GRM}. Then the order of $\gamma = q = m$.
\end{prop}
\begin{proof}
Since $G$ is type $III$ group, therefore any prime divisor of the order of $G$
is greater than or equal to $7$. Therefore if $q$ is not equal to $m$, then 
\begin{equation*}
q \leq \frac{m}{7}.
\end{equation*}
Using Theorem~\ref{GRM} we have the following equality for the density of 
the set $A$.
\begin{equation*}
\alpha(A) = \frac{1}{q} \sum_{i \in \Z/q\Z} \alpha_i =
   \frac{1}{q}\sum_{i \in I_q} \alpha_i
\end{equation*}
Now from  Lemma~\ref{2l}, it follows that
\begin{equation*}
\alpha(A) =  \frac{1}{q} \sum_{ i \in H_q \cup M_q} (\alpha_i + \alpha_{2i})
\leq  \frac{1}{q} \sum_{ i \in H_q \cup M_q} 1. 
\end{equation*}
Since the cardinality of the set $H_q \cup M_q$ is equal to $\frac{q - 1}{3}$,
it follows that 
\begin{equation*}
\alpha(A) \leq  \frac{1}{q}\frac{q - 1}{3} \leq \frac{1}{3} - \frac{7}{3m}.
\end{equation*}
But the last inequality above is contrary to assumption that  
\begin{equation*}
\alpha(A) > \mu(G) - \frac{5}{42m} > \frac{1}{3} -  
\frac{7}{3m}.
\end{equation*}
 Hence the lemma follows.
\end{proof} 
\section{Element with large fiber }\label{largefiber}
Given a set $A \subset G$ such that $\alpha(A) \geq \mu(G) - min( \eta, 
\frac{5}{42m})$, from Theorem~\ref{GRM} and Proposition~\ref{ord} it follows 
that $\gamma$ is a surjective homomorphism from $G$ to $\Z/m\Z$ such that
$A$ is a subset of the set $\gamma^{-1}(I)$,
where $\gamma$ is the special direction of the set $A$ and
 $I$ is a subset of $\Z/m\Z$ as defined in section~\ref{sintro}.
Then we define
$L \subset \Z/m\Z$ as follows.
\begin{equation*} 
L = \{ i \in \Z/m\Z:  \alpha(g, i) > \frac{1}{2} \}
\end{equation*}
We say that the fiber of an element 
$i \in \Z/m\Z$ is large if $i$ belong to the set $L(g, A)$. 
It is clear that $L$ is a subset of the set $I$
In this section we shall show that the set
$L$ is a sum-free subset 
 $\Z/m\Z$ and the cardinality of the set $L$ is $2k$, where
$k$ is equal to $\frac{m - 1}{6}$. 
\\
\\ 
The fact that the set $L$ is sum-free is a consequence of the following
folklore in additive number theory. We give a proof of the following lemma
for the sake of completeness.
\begin{lem}\label{llf} (folklore) Let $C$ and $B$ are subsets of a finite 
abelian group $K$ such that $min(|C|,|B|) > \frac{1}{2} |K|$. Then $C + B = K$.
\end{lem}
\begin{proof} Suppose there exist $x \in K$ such that $x$ does not
belong to $C + B$. This is clearly equivalent to the fact that $C \cap
(x - B) = \phi$. But this means that $|K| > |C| + |x - B| > |K|$. This
is not possible. Hence the lemma is true.
\end{proof}
 
\begin{lem}\label{ll0} For any two elements $i, j \in L$, we have
$\alpha_{i + j} = \alpha_{i - j} = 0$. In particular the set $L$ is sum-free.
\end{lem}
\begin{proof} The fact that $\gamma$ is a homomorphism implies that
the set $A_i + A_j$ is a subset of the set $ S_{i +j}$. Take any $x \in
S_i,  y \in S_j$.     Let $C = A_i - x$ and $B = A_j -
y$ so that we have the sets $C$ and $B$ are subsets of group  $S$. Then
applying Lemma~\ref{llf} it follows that $C + B = S  $. Therefore we have
$A_i + A_j = S_{i +j}$. The fact that $A$ is sum-free  implies that
$A_{i + j} \cap (A_i + A_j) = \phi$. Since we have shown that
the set $A_i + A_j = S_{i + j}$, it follows that the set $A_{i + j}= \phi$.
In other words it follows that  $\alpha_{i + j} =
0$. From similar arguments it also follows that $\alpha_{i - j} = 0$.
\end{proof}
Now we shall show that the cardinality of the set $L$ is equal to $2k$. For
this we require the following Lemma.
\begin{lem}\label{CHINEQ} Let $A, G,$ be as in theorem \ref{BGCH}. Let $m$
be the exponent of the group $G$ and $I, H, T, M$ are the subsets
of $\Z/m\Z$ as defined in section~\ref{sintro}. Let $m = 6k +1$.
Let $g$ be a surjective homomorphism $g: G \to \Z/m\Z$ such that the following 
holds.
\begin{equation*}
A \subset g^{-1} (I).
\end{equation*}
Then we have the following inequality
\begin{equation} \label{alb}
 \alpha(g, i) + \alpha(g, 2i) \geq   1 - m(\mu(G) - \alpha(A) ), 
\quad \forall i \in H \cup T.
\end{equation}
\end{lem}
\begin{proof} 
From the assumption on the homomorphism $g$ it follows that the density of the 
set $A$ satisfy the following equality
\begin{equation*}
\alpha(A) = \frac{1}{m}\sum_{i \in I} \alpha(g, i)
\end{equation*}
Now from Lemma~\ref{2l}, it follows that the set $I$ is divided into $2k$
disjoint pairs of the form $(i, 2i)$ where $i$ belongs to the set $H \cup T$.
Therefore it follows that
\begin{equation*}
\alpha(A) = \frac{1}{m}\sum_{i \in I} \alpha(g, i) = 
\frac{1}{m} \sum_{ i \in H \cup T} (\alpha(g, i) +  \alpha(g, 2i) ). 
\end{equation*}
Now using Lemma~\ref{alpha1} it follows that for any $i_0 \in H \cup T$ the 
following inequality holds
\begin{equation*}
m\alpha(A) \leq 2k - 1  + \alpha(g, i_0) + \alpha(g, 2i_0).
\end{equation*}
From this the required inequality follows for any $i_0$ belonging to
the set $H \cup T$ after observing that $\mu(G) = \frac{2k}{m}$. 
\end{proof} 
We need 
the following well known theorem due to Kneser.
\begin{thm} (Kneser) \label{knes} Let $C, B $ are subsets of a finite
abelian group $K$ such  that $|C + B | < |C| +  |B|$. Let $ F = H(C + B) =
 \{ g
\in G: g + C + B = C + B\}$ be the stabiliser of the set $C + B$. Then
the following holds
\begin{equation*} |C + B| = |C + F| + |B + F| - |F|.
\end{equation*} In particular the set $F$ is a nontrivial subgroup of
$K$ and 
\begin{equation}\label{stablb}
|F| \geq |C| +|B| - |C + B| .
\end{equation}
\end{thm}
\noindent For the proof of this theorem one may see \cite{Nath}. 
\begin{lem}\label{1s2b} Let $A, G,$ be as in theorem \ref{BGCH},
$\gamma$ be as provided by Theorem~\ref{GRM} and $H, T, M$ be as
defined earlier. Then the following holds
\begin{enumerate}
\item For any  $i \in H \cup T$ if $\alpha_i \leq
\frac{1}{2}$ then  $\alpha_{2i} > \frac{1}{2}$. 
\item The cardinality of the set $L$ is equal to $2k$.
\end{enumerate}
\end{lem}
\begin{proof} 
\begin{enumerate}
\item Suppose the claim is not true. Then there exist $i_0 \in H \cup T$ such 
that
\begin{eqnarray}
 \alpha_{i_0} &\leq& \frac{1}{2}\\
 \alpha_{2i_0} &\leq& \frac{1}{2}.  
\end{eqnarray}
Then from~\eqref{alb} it follows that
\begin{eqnarray}
 \alpha_{i_0} & > & \frac{1}{2} - \frac{5}{42} \label{ai} \\
 \alpha_{2i_0} & > & \frac{1}{2} - \frac{5}{42}. \label{a2i}  
\end{eqnarray}
Take any $x \in S_{i_0}$ and consider the set $A_{i_0} - x$. Then 
invoking~\eqref{disj} and using~\eqref{a2i}, it follows that
\begin{equation} \label{aiai}
|(A_{i_0} - x) + (A_{i_0} - x) | = |A_{i_0}  + A_{i_0} | \leq |S| - |A_{2i_0}|
< (\frac{1}{2} + \frac{5}{42}) |S|
\end{equation}
Therefore using~\eqref{ai} it follows that
\begin{equation} \label{ailb}
2|A_{i_0} - x| = 2|A_{i_0}| > 2(\frac{1}{2} - \frac{5}{42}) |S| > 
|(A_{i_0} - x) + (A_{i_0} - x) |
\end{equation} 
Let $F$ denote the stabiliser of the set $(A_{i_0} - x) + (A_{i_0} - x)$.
We can apply  Theorem~\ref{knes} with $C = B = A_{i_0} - x$ and
using~\eqref{stablb},~\eqref{aiai},~\eqref{ailb}  
we have the following inequality
\begin{equation}
|F| > \left(\frac{1}{2} - \frac{15}{42} \right)|S| = \frac{1}{7}|S|.
\end{equation}
Therefore the cardinality of the group $S/F$ is strictly less than $7$.
But since $S$ is a group of type $III$, the group $S/F$ is also
of type $III$. 
Hence it follows that $S = F$. Therefore the stabiliser of the set
$A_{i_0} - x$ is equal to the group $S$. This implies that the set 
$A_{i_0} = S_{i_0}$. This is in contradiction to the assumption that 
$\alpha_{i_0} \leq \frac{1}{2}$. Hence the claim follows.
\item The set $I$ is divided into $2k$ disjoint pairs of the form
$(i, 2i)$ with $i \in H \cup T$. From $(I)$ it follows that at least one 
element of any such pair belongs to the set $L$. The claim follows 
since we have shown that the set $L$ is sum-free and is a subset of the
 set $I$.
\end{enumerate}
\end{proof}
\noindent   
From Lemma~\ref{ll0} and Lemma~\ref{1s2b} the following proposition follows

\begin{prop} \label{Lf}
\begin{enumerate}
\item The set $L$ is a sum-free subset of $\Z/m\Z$ of cardinality $2k$. The 
set $L$ is a subset of the set $I$.
\item For any two elements $i, j \in L$, we have $\alpha_{i + j} = 
\alpha_{i - j} = 0$.
\end{enumerate}
\end{prop} 
\section{Sum-free subset of cyclic group }\label{scyclic}
Let the group $\Z/m\Z$ be of type $III$ group and $m = 6k + 1$.
In this section we shall characterise all the sets $E \subset \Z/m\Z$
such that the set $E$ is sum-free and $|E| = 2k$. From Theorem~\ref{GRM}
it is sufficient to characterise those sets $E$ which are subset of the set
$I$. 
\begin{lem} \label{Edisj}
Let $E \subset \Z/m\Z$ be a sum-free set. Let the group $\Z/m\Z$ is of 
type $III$ and the cardinality of the set $E$ is $2k$, where $k$ is equal to  
$\frac{m -1}{6}$. Let $H, T, M, I$ are subsets of $\Z/m\Z$ as defined in section~\ref{sintro}
and the set $E$ is a subset of the set $I$.  Then for any 
element $y$ belonging to the set $M$ exactly one of the element of the pair
$(\frac{y}{2}, y)$ belongs to the set $E$.
\end{lem}
\begin{proof}
This is straightforward from the fact that the set
$I$ is divided into $2k$ disjoint pairs of the form $(\frac{y}{2}, y)$ with
$y$ belonging to the set $M$ and the assumption that the set $E$ is a sum-free
set and  is a subset of the set $I$.
\end{proof}
We have the natural
projection $p$ from the set of integers to $\Z/m\Z$.
\begin{equation*} p: \Z \to \Z/m\Z
\end{equation*}
Since the map $p$ restricted to the set $\{0, 1, 2, \cdots, m -2, m -1\} 
\subset \Z$ is a bijection to the group $\Z/m\Z$ ( as a map of the sets)
we can define 
\begin{equation*} p^{-1}:  \Z/m\Z \to \{0, 1, 2, \cdots m -1 \}
\end{equation*} in an obvious way.\\
The following lemma
is straightforward to check.
\begin{lem} \label{propI}
\begin{enumerate}
\item The set $H$ is equal to the set $-T$. 
\item The set $M$ is equal to the set $-M$.
\item The set $I$ is equal to the set $-I$.
\item For any set $B \subset \Z/m\Z$ the set $B \cap T$ is same as the set
$-( (-B) \cap H)$. Also the set 
$p^{-1} (B \cap T) = m - p^{-1}( (-B) \cap H)$.
\item The set $H + H$ as well as the set $T + T$ are subsets of the set $M$.   
\item Given any even element $y$ belonging to the set $p^{-1}(M)$ the 
element $\frac{p(y)}{2}$ belong to the set $H$. Also the element
$p^{-1}(\frac{p(y)}{2})$ is equal to the element $\frac{y}{2}$.
\item Given any odd element $y$ belonging to the set $p^{-1}(M)$ the 
element $\frac{p(y)}{2}$ belong to the set $T$. Also the element
$p^{-1}(\frac{p(y)}{2})$ is equal to the element $\frac{y + 6k + 1}{2}$.
\item Given any two elements $x, y $ belonging to the set $\in p^{-1}(H)$
 which are of same parity, the element 
$p^{-1}(\frac{p(x + y)}{2}) = \frac{x + y}{2} $ and the element 
$\frac{x + y}{2}$ belongs to the set $p^{-1}(H)$. Moreover if $x$ is strictly less 
than $y$ then the following inequality holds.
\begin{equation*}
x < \frac{ x + y}{2} < y
\end{equation*}
\end{enumerate}
\end{lem}
\begin{lem} \label{mid}
Let $E$ be a set as above. Then the following holds.
\begin{enumerate}
\item Given any two elements $x, y$ which belong to
the set $p^{-1}(E \cap H)$ and are of same parity, the element 
$\frac{x + y}{2}$ belong to the set $p^{-1}(E \cap H)$.
\item Given any two elements $x, y$ which belong to
the set $p^{-1}(E \cap H)$ and are of different parity, the element 
$\frac{x + y + 6k + 1}{2}$ belong to the set $p^{-1}(E \cap T)$.
\item Given an element $x$ belonging to the set $p^{-1}(E \cap H)$ and
an element $y$ belonging to the set $p^{-1}(E \cap T)$, the element
$\frac{p(x) -p( y)}{2}$ belong to the set $p^{-1}(E \cap (H \cup T))$. 
\item  Any two consecutive element
of the set $p^{-1}(E \cap H)$ (or of the set $p^{-1}(E \cap T)$ )are of 
different parity.
\end{enumerate}
\end{lem}
\begin{proof}
Since the set $E$ is sum-free, it follows that given
 any two elements $x, y$ belonging to the set $p^{-1}(E)$, neither the
element $p(x)~ + ~p(y)$ nor the element  $p(x)~ -~ p(y)$ belong to the
set $E$. Using this we prove all the claims.
\begin{enumerate}
\item Under the assumption, the element $p(x)~ + ~p(y)$ belong to the set $M$. 
From Lemma~\ref{Edisj} it follows that the element $\frac{p(x)~ + ~p(y)}{2}$
belong to the set $E$.
Also the element $p^{-1}(p(x) + p(y))$ is equal to $x + y$ and is even.
Therefore invoking Lemma~\ref{propI} it follows that 
the element $p^{-1}\left(\frac{p(x)~ + ~p(y)}{2}\right)$ is equal
to $\frac{x + y}{2}$ and belong to the set $p^{-1}( H)$. Hence the claim 
follows.
\item Under the assumption, the element  $p(x)~ + ~p(y)$ belong to the set
  $M$. From Lemma~\ref{Edisj} it follows that the element 
$\frac{p(x)~ + ~p(y)}{2}$
belong to the set $E$. In this case the element $p^{-1}p(x + y)$ is equal to 
the element  $x + y$ and is odd. Therefore invoking Lemma~\ref{propI} it
 follows that the element
$p^{-1}\left(\frac{p(x)~ + ~p(y)}{2}\right)$ is equal to the element 
$ \frac{ x + y + 6k + 1}{2}$ and belongs to the set $p^{-1}( T)$. Hence 
the claim follows.
\item Under the assumption, the element $p(x)~ - ~p(y)$ belong to the set $M$.
Therefore the claim follows invoking Lemma~\ref{Edisj}.
\item Let the set $p^{-1}(E \cap H) = \{ x_1 < x_2 < \cdots < x_h \} $.
Suppose 
there exist $ 1 \leq i_0 \leq h - 1$ such that the element $x_{i_0}$ and
the element $x_{i_0 + 1}$ have same parity. Then from $(I)$ it follows that
the element $\frac{x_{i_0} + x_{i_0 + 1}}{2}$ belong to the set 
$p^{-1}(E \cap H)$. From Lemma~\ref{propI} the following inequality also
 follows. 
\begin{equation*}
x_{i_o} < \frac{ x_{i_o} + x_{i_o + 1}}{2} < x_{i_o + 1} 
\end{equation*}
But this contradicts the fact that the elements 
$x_{i_o}$ and $x_{i_o + 1}$ are consecutive 
elements of the set $p^{-1}(E \cap H)$. Therefore the claim follows for the
set $p^{-1}(E \cap H)$. Replacing the set $E$ by the set $-E$, it follows that
the any two consecutive element of the set 
$p^{-1}\left((-E) \cap H \right)$ are also of different
parity.
Noticing that  the set
$p^{-1}(E \cap T) = m - p^{-1}\left((-E) \cap H \right)$, the claim follows 
for the set
$p^{-1}(E \cap T)$ also. 
\end{enumerate}
\end{proof}
\begin{prop} \label{ap}
Let $E$ be a set as above then the following holds.
\begin{enumerate}
\item The set $p^{-1} (E \cap H)$ as well as the set $p^{-1}(E \cap T)$ 
is an arithmetic progression with an odd common difference. 
\item The following inequality holds.
\begin{equation}
|E \cap T| -1 \leq |E \cap H| \leq |E \cap T| + 1
\end{equation}
\end{enumerate}
\end{prop}
\begin{proof}
Let the set $p^{-1}(E \cap H) = \{ x_1 < x_2 < \cdots < x_h \} $. 
\begin{enumerate}
\item In case the cardinality of the set $E \cap H$ is less than or equal to 
$2$, the claim is trivial for the set $p^{-1} (E \cap H)$. Otherwise
for any $ 1 \leq i \leq h -2$, consider the elements 
$x_{i}, x_{i + 1}, x_{i + 2}$, then from Lemma~\ref{mid} it follows that
the parity of elements $x_{i}$ and $x_{i +1}$ are different. For the same
reason the the parity of elements $x_{i + 1}$ and $x_{i +2}$ are different.
Therefore the parity of elements $x_{i}$ and $x_{i +2}$ are same. Therefore 
from Lemma~\ref{mid} it follows that the element $\frac{x_i + x_{i +2} }{2}$ 
belong to the set $p^{-1} (E \cap H)$. But from Lemma~\ref{propI} the
following inequality follows.
\begin{equation*}
x_i < \frac{x_i + x_{i +2} }{2} < x_{i + 1}
\end{equation*}
Hence it follows that for for any $ 1 \leq i \leq h -2$
\begin{equation*}
\frac{x_i + x_{i +2} }{2} = x_{i + 1}.
\end{equation*} 
This is equivalent to the fact that the set $p^{-1} (E \cap H)$ is an
arithmetic progression. It also follows that the common difference is odd.
The claim  for the set 
$p^{-1} (E \cap T)$ follows by replacing the set $E$ by the set $-E$.
\item From Lemma~\ref{mid} it follows that the set
\begin{equation} \label{ETSub}
\{ \frac{x_1 + x_2 + 6k + 1}{2} < \frac{x_2 + x_3 + 6k + 1}{2} < \cdots <
\frac{x_{h -1} + x_h + 6k + 1}{2} \}
\subset p^{-1}(E \cap T).
\end{equation}
 Therefore it follows that
\begin{equation*}
|E \cap T| \geq |E \cap H| - 1 .
\end{equation*} 
Replacing the set $E$ by the set $-E$, it also follows that
\begin{equation*}
|E \cap H| \geq |E \cap T| - 1 .
\end{equation*} 
Hence the claim follows.
\end{enumerate}
\end{proof}
\subsection{$max(|E \cap H|, |E \cap T|) \geq 2$ }
\begin{prop}\label{int}
Let $E$ be a set as above and $H, T, M$ as defined above. 
\begin{enumerate}
\item  Suppose the inequality $min(|E\cap H|, |E\cap T|) \geq 2$ is satisfied.
Then the set $p^{-1}(E \cap H)$ and the set $E \cap T$ are arithmetic 
progression
with same common difference $d(H, E) = d(T, E) = d(E) 
(\text{ say })$. 
\item  Suppose the inequality $min(|E\cap H|, |E\cap T|) \geq 2$ is satisfied.
The set $p^{-1}(E^c \cap M)$ is an arithmetic
progression
with common difference $d(E)$, where $d(E)$ is a positive integer given by 
$(I)$.
\item Suppose the cardinality of the set $p^{-1}(E \cap H)$ 
(resp. $p^{-1}(E \cap T)$ ) is equal to $2$ and 
the cardinality of the set $p^{-1}(E \cap H)$ 
(resp. $p^{-1}(E \cap T)$ ) is equal to $1$, then  
the set $p^{-1}(E^c \cap M)$ is  an arithmetic
progression
with common difference $d(H,E)$ (resp. $d(T, E)$ ) which is equal to $1$.
 \item Let the inequality $max(|E\cap H|, |E\cap T|) \geq 2$ is satisfied. 
Then 
the set $p^{-1}(E^c \cap M)$ is  an arithmetic
progression
with common difference equal to $1$.
\end{enumerate}
\end{prop}
\begin{proof}
\begin{enumerate}
\item We discuss the two cases.\\
{\bf Case 1:} $max ((|E\cap H|, |E\cap T|) \geq 3$.\\
Under the assumption, either the inequality $|E\cap H| \geq 3$ holds
or the inequality $|E\cap T| \geq 3$ holds (both the inequalities may also
hold). Since
 the claim holds for the set 
$E$ if and only if it holds for the set $-E$, it is sufficient to prove the
assertion under the assumption that   the inequality $|E\cap H| \geq 3$
holds.\\
Now from Proposition~\ref{ap} it follows that the the sets $p^{-1}(E \cap H)$ 
and  $E \cap T$ are arithmetic 
progression. Let the set
\begin{equation*} 
p^{-1}(E \cap H) = \{x_1 < x_1 + d(H, E) < x_2 + 2d(H, E) < \cdots < 
x_1 + (h -1) d(H, E) \}.
\end{equation*}
Then from~\eqref{ETSub} it follows that the  following set 
\begin{equation*}
\{x_1 +  \frac{ d(H,E) + 6K + 1}{2}, x_1 + \frac{ d(H,E) + 6K + 1}{2} + 
d(H, E)\}
\end{equation*}
 is a subset of the set $p^{-1}(E \cap T)$. From this the following inequality
follows immediately
\begin{equation*}
d(H, E) \leq d(T, E).
\end{equation*}
In the case the cardinality of the set $E \cap T$ is equal to $2$ it also 
follows that the set 
\begin{equation*} 
 \{x_1 +  \frac{ d(H,E) + 6K + 1}{2}, x_1 + \frac{ d(H,E) + 6K + 1}{2} + 
d(H, E)\}
\end{equation*}
is equal to the set $p^{-1}(E \cap T)$. Hence the claim follows in case
the cardinality of the set $E \cap T$ is equal to $2$.
Suppose the cardinality of the set $E \cap T$ is also greater than or
equal to $3$, then replacing the set $E$ by the set $-E$, it follows that
\begin{equation*}
d(H, -E) \leq d(T, -E).
\end{equation*} 
Since the numbers $d(H, -E)$ and $d(T, -E)$ are equal to the numbers 
$d(T, E)$ and $d(H, E)$ respectively, the claim follows.\\
\\
\noindent
{\bf Case 2:} $|E \cap H| = |E \cap T| = 2$. \\
Let the set $p^{-1}(E \cap H) = \{x, y\}$ and the set  
$p^{-1}(E \cap T) = \{z, w\}$. Then from Lemma~\ref{mid}, it follows that
the parity of the elements $x$ and $y$ are different and the element
$\frac{ x + y + 6k + 1}{2}$ belong to the set  
$p^{-1}(E \cap T)$. Since we are not assuming that $z < w$ we can assume
without any loss of generality that
the element 
\begin{equation*}
\frac{ x + y + 6k + 1}{2} = z.
\end{equation*}
For the similar reason the element 
 $\frac{ z + w - 6k + 1}{2}$ belong to the set  
$p^{-1}(E \cap T)$ and we can  assume
without any loss of generality that
the element 
\begin{equation*}
\frac{ z + w - 6k + 1}{2} = x.
\end{equation*}
Therefore it follows that
\begin{equation*}
w -z = x - y.
\end{equation*}
This proves the claim.
\item 
First we notice that the claim is true for the set $E$ if and only if it is
true for the set $-E$. This is because the sets $(-E)^c$ and
$M$ are equal to the  sets   $-(E)^c$ and $-M$ respectively. Therefore
the set $p^{-1}(E^c \cap M)$ is same as the set 
$ 6k + 1 - p^{-1}((-E)^c \cap M)$. Therefore replacing the set $E$ by the set
$-E$ if necessary we may assume that the following inequality holds.
\begin{equation*}
|E \cap H| \geq |E \cap T|
\end{equation*}
 Let the smallest member of 
$p^{-1}(E \cap H)$ be $x$ so that 
$p^{-1}(E \cap H)= \{x + jd(E):0\le j\le h-1\}$.
Using~\eqref{ETSub} it follows that the  set
\begin{equation}
\{a+ jd(E):0\le j\le h-2\}
~~~~~{\rm where }~~ a=x + 3k + \frac{d(E) + 1}{2},
\end{equation}
is a subset of  $p^{-1}(E \cap T)$ and its cardinality  is $h-1$. 
Since the cardinality of $p^{-1}(E \cap T)$ is at-most $h$ and 
since $p^{-1}(E \cap T)$ is in arithmetic progression with common difference
$d(E)$, it follows that it is contained in 
$\{a+ jd(E):-1\le j\le h-1\}$. On the other hand,  $p^{-1}(E^c \cap M)$ 
is the disjoint  union of
$\{2x: x \in p^{-1}(E \cap H) \}$ and $\{2x - 6k -1: x  
\in p^{-1}(E \cap T) \}$. Therefore, we have
$\{2x + jd(E): 0 \le j \le 2h -2\} \subset p^{-1}(E^c \cap M)
\subset \{2x + jd(E): -1 \le j \le 2h -1\}.$
From this the claim is immediate.
\item It is sufficient to prove the assertion in the case the cardinality of 
the sets $p^{-1}(E \cap H)$ and $p^{-1}(E\cap T)$ are equal to $2$ and $1$
respectively. In the other case then the assertion follows by replacing the
set $E$ by the set $-E$. Suppose the set $p^{-1}(E \cap H)$ is $\{x, y\}$
and $p^{-1}(E \cap T)$ is $\{z\}$. Then it follows that the element
 $2z - 6k - 1$ is equal to the element $ x + y$. It also follows that the
set $p^{-1}(E^c \cap M)$ is $\{x, x + y, 2y\}$. Hence the assertion follows.
\item Replacing the set $E$ by the set $-E$ we may assume $|E \cap H| \geq 2$.
 From $(I), (II)$ and $(III)$ we know 
that the set $p^{-1}(E^c \cap M)$ is an arithmetic progression with common
difference equal to $d(H. E)$ which is an odd positive integer.
Suppose the assertion is not true. Then it means that $d(H, E) \ge 3$.
Then the smallest element of the set  $p^{-1}(E \cap H)$ (let say $x$)
is less than or equal
to $2k - 3$. This implies that $k$ is at least $3$.
It also follows that the cardinality of the set 
$\{2k + 1, 2k + 2, 2k + 3\} \cap (E \cap M) \ge 2$. But since $ x \le 2k -3$, 
the set $\{ x + 2k + 1, x +  2k + 2, x + 2k + 3\}$ is a subset of the set 
$p^{-1}(M)$
and contains at least $2$ elements of the set $p^{-1}(E^c \cap M)$. But this
contradicts the fact that the set   $p^{-1}(E^c \cap M)$ is in arithmetic 
progression of common difference $d(H, E)$ which is at least $3$. Hence the
claim follows.
\end{enumerate}
\end{proof}
\begin{defi}
 We shall say that a set $B \subset Z$ is an interval if
 either it is an arithmetic progression with common difference equal to $1$ or 
the cardinality of the set $B$ is equal to $1$.
\end{defi}
\begin{prop} \label{HUT}
Let $E$ be a set as above.
 Suppose ${\text max}(|E \cap H|, |E \cap T| ) \geq 2$ then
$ E = H \cup T$.
\end{prop}
\begin{proof}
From Lemma~\ref{Edisj}, the conclusion of proposition is 
equivalent of the assertion that $p^{-1}(E^c \cap M) = p^{-1}(M)$. Let
$y$ be the smallest member of $p^{-1}(E^c \cap M)$ and it's cardinality is
$s$. Then using Proposition~\ref{int} it follows $p^{-1}(E^c \cap M)$
is equal to $\{y + j: 0 \le j \le s -1 \}$.  Suppose the claim is not true.
Therefore 
 at~least one of the element of the set $\{2k +1, 4k\}$  belong to
the set $p^{-1}(E)$.
 The assumption of the proposition is satisfied
for the set $-E$ as well and the conclusion is true for the set $E$ if and only
if it is true for the set $-E$. Therefore replacing
the set $E$ by the set $-E$ if necessary
we may assume that the element $2k + 1$ belong
to the set $p^{-1}(E)$. Using this we  prove the following.\\
{\bf Claim:} $p^{-1}(E \cap H) = p^{-1}(H)$.\\
Let
$x$ be the smallest member of $p^{-1}(E^c \cap M)$ and it's cardinality is
$h$. Then using Proposition~\ref{int} it follows $p^{-1}(E \cap H)$
is equal to $\{x + j: 0 \le j \le h -1 \}$. Therefore it is sufficient to show 
that $x + h - 1$ is $2k$ and $x$ is $k + 1$.\\
\\
\noindent
{\bf Claim} The largest element of $p^{-1}(E \cap H)$ is $2k$.\\
\\
\noindent
Suppose $x + h - 1$ is not equal to $2k$. 
Since the set
$E$ is sum-free, it follows that $\{2x + j: 0 \le j \le 2h -2\}$ is a subset
of  $p^{-1}(E^c \cap M)$.  As $x + h -1 $ is not $2k$ therefore $2x + 2h$
belong to  $p^{-1}(E \cap M)$. Therefore $y + s -1 \le 2x + 2h -1$. 
Since $2k + 1$ belong to the set $p^{-1}(E)$ it also follows
that $ x + h - 1 + 2k + 1$ belong to the set $p^{-1}(E^c)$. Since 
$ x + h - 1 < 2K$, we have $x + h -1 + 2k + 1$ belong 
to the set $p^{-1}(E^c \cap M)$. Therefore it follows that
$ x + h -1 + 2k + 1 \le 2x + 2h - 1$ which implies $ x + h -1 \ge 2k$. Since
$ x + h -1$  belong to the set $p^{-1}(H)$ it follows that $x + h -1 = 2k$, 
contradictory to the assumption $   x + h -1 \ne 2k$. Therefore $x + h -1$ is
$2k$.\\
\\
\noindent
{\bf Claim} The smallest element of $p^{-1}(E \cap H)$ is $k + 1$.\\
\\
\noindent
Now we will show that $x$ is $k + 1$. Since we have shown that
$p^{-1}(E \cap H) = \{x + j: 0 \le j \le 2k - x\}$, using Lemma~\ref{mid} it 
follows that $\{x + j : 3k + 1 \le j \le 5k - x\}$ is a subset of 
$p^{-1}(E \cap T)$ and it's cardinality is $h -1 = 2k -x$.  The set
$p^{-1}(E \cap T)$ is an interval and it's cardinality is at~most $h + 1$.
Therefore we have that $p^{-1}(E \cap T)$ is a subset of
 $\{x + j : 3k -1 \le j \le 5k - x\}$. But in case $\{ x + 3k -1, x + 3k +
1\}$
is a subset of $p^{-1}(E \cap T)$, from Lemma~\ref{mid} it follows that
$x - 1$ belong to $p^{-1}(E \cap T)$. This contradicts that $x$ is the
smallest
element of $p^{-1}(E \cap H)$. Therefore $p^{-1}(E \cap T)$ is a subset of
 $\{x + j : 3k  \le j \le 5k - x\}$. 
Therefore the least element of$p^{-1}(E \cap T)$ is either 
$ x + 3k$ or $x + 3k +1$.\\
\\
{\bf Case} The  least element of $p^{-1}(E \cap T)$ is  $x + 3k$\\
\\
\noindent
In this case the element $x + 3k - 2k = x + k$ does not belong to the 
set $p^{-1}(E)$. Since $x$ belong to $p^{-1}(H)$, the element $x + k$ belong
to $p^{-1}(M)$. Therefore in case $x + k$ is even, the element 
$\frac{x + k}{2}  $ belong to the set $p^{-1}(E \cap H)$ . But if $x + k$ 
is even then  $\frac{x + k}{2} < x  $. This contradicts  the assumption
that the element
$x$ is the least element of $p^{-1}(E \cap H)$. In case $x + k$ is odd, then
the element $\frac{x + k + 6k + 1}{2}  $ belong to the set $p^{-1}(E \cap T)$.
Now in case $x$ is not $k + 1$ then $ x \ge k+ 2$ and
 the inequality $\frac{x + k + 6k + 1}{2} < x + 3k $ is  satisfied.
This contradicts that   $ x + 3k$ is the least element of  $p^{-1}(E \cap T)$.
Therefore $x$ has to be $k +1$ in this case.
\\
\\
\noindent
{\bf Case} The  least element of $p^{-1}(E \cap T)$ is  $x + 3k + 1$.\\
\\
\noindent
In
 this case
the element $x + k + 1$  belong to the set $p^{-1}(E^c \cap M)$.
Then again if $x$ is not $k + 1$, this leads to a contradiction.\\
\\
\noindent
Therefore  $p^{-1}(E \cap H) = H$. The assumption of proposition can hold
only in case $k \ge 2$. In that case $k + 1$ and the element $k + 2$ belong to
$p{-1}(H) = p^{-1}(E \cap H)$. Therefore the element $2k + 3$ does not belong 
to
$p^{-1}(E)$ and invoking Lemma~\ref{Edisj} the element 
$p^{-1}p\left(\frac{2k + 3}{2}\right) = 4k + 2$
belong to $p^{-1}(E \cap T)$. But since by assumption $2k + 1$ belong to the
set $p^{-1}(E)$, this contradicts that the set $E$ is sum-free.
Therefore finally it follows that the set $p^{-1}(E^c \cap M) = M$ and hence
the proposition follows.
\end{proof}
\subsection{ $max(|E \cap H|, |E \cap T|) \le 1$}
Replacing the set $E$ by the set $-E$ if necessary we may assume that the 
following inequality holds.
\begin{equation*}
|E \cap H| \ge |E \cap T|
\end{equation*}
Then we have the following three possible cases.
\begin{enumerate}
\item The equality $|E \cap H| = |E \cap T| = 1$ holds.
\item The cardinality of the sets $E \cap H$ and $E \cap T $ are $1$ and
$0$ respectively.
\item The equality $|E \cap H| = |E \cap T| = 0$.
\end{enumerate}

\begin{prop} \label{11}
Let $E, H,T, M $ be as above. If $|E \cap H| = |E \cap T| = 1$, then the set
$p^{-1}(E) = \{2k\} \cup \{2k+2, 2k + 3, \cdots, 4k - 2, 4k -1\} \cup \{4k + 1 \}$.
\end{prop}
\begin{proof}
Let the set $(E\cap H) = \{x\}$ and the set $(E \cap T) = \{y\}$. 
From 
Proposition~\ref{Edisj} the
set $E^c \cap M = \{2x, 2y \}$ and the set $p^{-1}(E^c \cap M) = 
\{p^{-1}(2x), p^{-1}(2y)
\}$. 
We claim that
 \\
\\
 {\bf claim :} x = -y.\\
{\it proof of claim: } From Proposition~\ref{Edisj} it follows that the
element $\frac{ -y + x }{2}$ belongs to the set $E\cap (H \cup T)$. Now if
 $p^{-1}(-y)$ and $p^{-1}(x)$ have same parity then the element  
$p^{-1}(\frac{ -y + x }{2})$ belongs
 to the set $p^{-1}( E \cap H)$. 
But as the element $p^{-1}(x)$ is the only element belonging to 
the 
set $E \cap H$, in this case the claim follows. Otherwise the element 
$\frac{ -y + x }{2}$ belongs to the set $E \cap T$ and hence is equal to $y$.
Also then $p^{-1}(-x)$ and $p^{-1}(y)$ have different parity and the element 
$\frac{ -x + y }{2}$ is equal to $x$. But this implies after simple
 calculation
that $x = 9x$. As $m$ is odd this is not possible. Hence the claim follows.\\
\\
Next claim is \\
{\bf claim: } $p^{-1}(x) = 2k$ \\
{\it proof of claim: } Suppose not, then $p^{-1}(y)$ is also not equal to 
$4k +1$. Therefore the element $2k + 1$ belong to the set $p^{-1}(E)$.
This implies that the element $p^{-1}(x) + 2k +1$ belongs to the set 
$p^{-1}(E^c \cap M)$.
The element  $p^{-1}(x) + 2k +1$ also satisfy the inequality 
$p^{-1}(x) + 2k +1 > 2p^{-1}x$. 
Therefore the element $p^{-1}(x) + 2k +1 = p^{-1}(2y)$. Since $y = -x$, it
follows that $p^{-1}(2y) = 6k + 1 - 2p^{-1}(x)$. Therefore we have 
$3p^{-1}(x) = 4k$. This is possible only if $k$ is divisible by $3$. It is easy
to check that case $k =3$ is not possible.
So we may assume that $k$ is greater than  $3$ and is divisible by
$3$. As $k$ is strictly greater than  $3$ therefore 
$p^{-1}(x) = \frac{4k}{3} \neq k + 1$
and hence $2k + 2$ belong to the set $E$.
Also we have the inequality $\frac{4k}{3} \leq 2k -2$. Therefore
the element $p^{-1}(x) + 2k + 2$ belong to the set $E^c \cap M$. Therefore the
elements
$ p^{-1}(x) + 2k + 1$ as well as $p^{-1}(x) + 2k + 2$ belong to the set 
$p^{-1}(E^c \cap M)$ and
neither of these elements are same as $2p^{-1}(x)$. This implies that the 
cardinality 
of the set $p^{-1}(E^c \cap M)$ is greater than or equal to $3$. 
This is not possible.
Hence the claim follows.\\
Now the proposition follows immediately.
\end{proof}
\begin{prop} \label{10}
Let $E$ be as set as above, then following holds.
If $|E\cap H| = 1$ and $|E \cap T| = 0$ then we have the set 
$p^{-1}(E) = \{2k,2k+1, \cdots, 4k - 2, 4k - 1\}$.
\end{prop}
\begin{proof}
Suppose the set $\{x\}$ is $ p^{-1}(E \cap H)$
The claim is immediate from the assertion $x = \{2k\}$.
Suppose the assertion is not true. Since we have assumed $|E \cap T| = 0$,
using Lemma~\ref{Edisj} it follows that $2k + 1$ belong to the set $p^{-1}(E)$.
Therefore if  $x \ne \{2k\}$, then $x + 2k + 1 $ belong to the set $p^{-1}(M)$
and actually belong to 
$p^{-1}(E^c \cap M)$. But trivially $p^{-1}(E^c \cap M) = \{2x\}$ and the 
element $x + 2k + 1$ is not equal to $2x$. Hence there is a contradiction and
the claim follows.
\end{proof}
The following proposition is trivial.
\begin{prop} \label{00}
Let $E$ be a set as above.
In the case $|E \cap H| = |E \cap M| = 0$, then $E = M$ 
\end{prop}
Therefore from Proposition~\ref{HUT},~\ref{11},~\ref{10} the proof of 
Theorem~\ref{maxch} in case $G$ is cyclic follows. That is the following 
result 
follows.
\begin{thm} \label{cyclic}
Let $G$ be a cyclic abelian group of type $III$. That is $G = \Z/m\Z$. Let
$k = \frac{m -1}{6}$. Let $E$ be a sum-free subset of $G$ of density $\mu(G)$.
Then there exist an automorphism $f: G \to G$ such that one of 
the following holds.
\begin{enumerate}
\item The set $E = f^{-1}(M)$.
\item 
The set $E$  is equal to $f^{-1}p( \{2k\} \cup \{2k + j: 2 \le j \le
 2k -1\}  \cup \{4k + 1\})$.
\item The set $E$ or the set $-E$ is equal to $f^{-1}p( \{2k\} 
\cup \{2k + j: 1 \le j \le
 2k -1 \}  )$.
\end{enumerate}
\end{thm}
\begin{proof}
From Theorem~\ref{GRM} and Proposition~\ref{ord}, it follows that there exist
automorphism $g: G \to G$ such that the set $E$ is a subset of the set
$g^{-1}(I)$. Then in case $max(|g(E) \cap H|, |g(E) \cap T|) \geq 2$ 
Proposition~\ref{HUT} we have $g(E) = H \cup T$. Therefore taking $f = 2g$ it 
follows $E = f^{-1}(M)$. In the other cases taking $f = g$ the claim follows
from Proposition~\ref{11},~\ref{10},~\ref{00}. 
\end{proof}
\section{Sum-free set of general abelian group of type III}\label{sgeneral}
\begin{proof}{\it  of Theorem~\ref{BGCH}:} 
Let $\gamma$ be a special direction
of the set $A$.  Then from Proposition~\ref{ord} it follows that order 
of $\gamma$ is equal to $m$.
We have also from Theorem~\ref{GRM}  it
follows that $A$ is a subset of $\gamma^{-1}(I)$. Let $L \subset \Z/m\Z$ as
defined in section~\ref{largefiber}. Then from Proposition~\ref{Lf} we have 
that $L$ is a sum-free subset of $\Z/m\Z$ and it's cardinality is $2k$.
Then from propositions~\ref{HUT},~\ref{11},~\ref{10},~\ref{00} it follows 
$p^{-1}(L)$ is one of the following set.\\
{\bf Case 1}: $p^{-1}(L) = p^{-1}(H \cup T)$. \\
In this case  we easily check that given any element $i \in \Z/m\Z$ such that
$i$ does not belong to the set $L$, there exist $x, y \in L$ such that 
either $i = x + y$ or $i = x -y$. Therefor invoking Proposition~\ref{Lf}, it
follows that for any $i \in \Z/m\Z$ which does not belong to
$L$ , the number $\alpha_i = 0$. Therefore $A$ is a subset of
$\gamma^{-1}(I)$. Hence $f = 2\gamma$ is a surjective homomorphism from
$G$ to $\Z/m\Z$ and moreover
\[
A \subset f^{-1}\{M\}.
\]
{\bf Case 2}:  $p^{-1}(L) = \{2k\} \cup \{2k + j: 2 \le j \le
2k -1 \} \cup \{4k + 1 \}$.\\
In this case again given any element $i \in \Z/m\Z$ such that
$i$ does not belong to the set $L$, there exist $x, y \in L$ such that 
either $i = x + y$ or $i = x -y$. Therefore taking $f = \gamma$ we have
\[
A \subset f^{-1}p(2k) \cup f^{-1}p(\{2k + j: 2 \le j \le
2k -1 \})  \cup f^{-1}p(4k + 1).
\]
\\
{\bf Case 3}:Either $p^{-1}(L)$ or $p^{-1}(-L)$ is
$\{2k\} \cup \{2k+j: 1 \le j \le 2k -1\}$.\\
In both these cases, given any element $i \in \Z/m\Z$ such that
$i$ does not belong to the set $L$, there exist $x, y \in L$ such that 
either $i = x + y$ or $i = x -y$. Therefore taking $f = \gamma$ the claims
follows.\\
\[
A \subset f^{-1}p(2k) \cup f^{-1}p(\{2k + j: 2 \le j \le
2k -1 \})  \cup f^{-1}p(4k + 1).
\]
{\bf Case 4}:$p^{-1}(L) = p^{-1}(M)$.\\
In this case given any $i \in \Z/m\Z$ which does not belong to 
$\{2k + j: 0 \le \le 2k + 1\}$ there exist  $x, y \in L$ such that 
either $i = x + y$ or $i = x -y$. Therefore taking $f = \gamma$ we have
\[
A \subset  f^{-1}p(\{2k + j: 0 \le j \le
2k + 1 \})
\]
Therefore in all the cases there exist a homomorphism $f: G \to \Z/m\Z$ such 
that 
\[
A \subset  f^{-1}p(\{2k + j: 0 \le j \le
2k + 1 \}).
\]
Now the claims $(I), (II), (III)$ follows invoking Lemma~\ref{CHINEQ}.
The claim $(IV)$ follows observing that the set 
$A \cup f^{-1}p(\{2k + j: 0 \le j \le
2k + 1 \}$ is also a sum-free set.
\end{proof}
Now we prove Theorem~\ref{maxch}.
\begin{proof} {\it of Theorem~\ref{maxch}}: Let $f$ be a surjective homomorphism
from $G$ to $\Z/m\Z$ given by Theorem~\ref{BGCH}. 
Since $m$ is the exponent of group $G$ and $f$ is a surjective homomorphism to
$\Z/m\Z$ there exist a subgroup $C$ of $G$ such that $G = S \oplus C$ where
$S$ is a kernel of $f$. Therefore $f$ restricted to $C$ is an isomorphism from 
$C$ to $\Z/m\Z$. We denote this restriction by $g$.
Since $\alpha(A) = \mu(G)$ it follows
from Theorem~\ref{BGCH} we have
\[
f^{-1}p(\{2k + j: 0 \le j \le
2k + 1 \}) \subset A.
\] 
Also the following equalities hold.
\begin{eqnarray}
|A(f, 2k)| + |A(f, 4k)| &=& \frac{n}{m} \label{2kcos}\\
|A(f, 4k +1)| + |A( f, 2k + 1)| &=& \frac{n}{m} 
\label{4kcos}
\end{eqnarray}  
Now  $A(f, 2k) + A(f, 2k)$ is a subset of $f^{-1}\{4k\}$ and it is
disjoint from  $A(f, 4k)$. Therefore the following inequality follows.
\begin{equation*}
|A(f, 4k)| \le \frac{n}{m} -  |A(f, 2k) + A(f, 2k)| = |A(f, 4k)| + |A(f, 2k)|
-  |A(f, 2k) + A(f, 2k)|
\end{equation*}
Hence we have
\begin{equation}\label{2keq}
|A(f, 2k) + A(f, 2k)| = |A(f, 2k)|. 
\end{equation}
For any $i \in \Z/m\Z$ there exist $X_i \subset S$ such that $A(f, i) =
X_i \oplus g^{-1}\{i\}$. 
Then from~\eqref{2keq} it follows
that $|X_{2k} + X_{2k}| = |X_{2k}|$. Therefore either $X_{2k} = \phi$ or 
there exist $J_1$ a 
subgroup of $S$ and an element $b_1 \in S$ such that $X_{2k} = J_1 + b_1$.
 Similar arguments implies that either  
$X_{4k + 1} = \phi$ or there exist $J_2$ a 
subgroup of $S$ and an element $b_2 \in S$ such that $X_{4k + 1} = J_2 + b_1$.
  Then there are three possibilities.\\
\\
{\bf Case 1}: Both the sets $X_{2k}$ as well as the set $X_{4k + 1}$
are empty sets.\\
In this case from~\eqref{2kcos} and~\eqref{4kcos} it follows
$X_{2k + 1} = S$ and $X_{4k} = S$. Hence 
$A$ is $f^{-1}(M)$.\\
\\
{\bf Case 2}: Exactly one of the sets $X_{2k}$ and $X_{4k + 1}$ is 
an empty 
set.\\
Replacing the set $A$ by $-A$ if necessary we may assume that 
$X_{4k + 1}$ is an empty set. Since the set $A$ is sum-free it follows that
$X_{4k}$ is a subset of $(J_1  + 2b_1)^c$.  From~\eqref{2kcos} it follows 
trivially
that $X_{4k}$is  $(J_1  + 2b_1)^c$.\\
{\bf Case 3}: Both the sets $X_{2k}$ and $X_{4k + 1}$ are not empty 
sets.\\
Then arguing as in case 2 it follows that $X_{4k}$ is  $(J_1  + 2b_1)^c$
and $X_{4k + 1}$ is  $(J_2  + 2b_2)^c$.
 The assumption that $A$ is sum-free implies that
 $X_{4k + 1}$
is a subset of $(X_{2k} + X_{2k + 1})^c$. This means
\[
J_2 + b_2 \subset (J_1 + b_1 + (J_2 + 2b_2)^c)^c.
\] 
This implies 
\[
(J_2 + b_2)^c = J_1 + b_1 + (J_2 + 2b_2)^c.
\]
Therefore we have
\begin{equation}\label{j1sub}
J_1 + b_1 \subset J_2 - b_2.
\end{equation}
Since $X_2k$ is a subset of $(X_{2k} + X_{2k + 1})^c$, same
arguments implies 
that
\begin{equation}\label{j2sub}
J_2 + b_2 \subset J_1 - b_1.
\end{equation}
From~\eqref{j1sub} and~\eqref{j2sub} we have $J_1 + b_1 = J_2 - b_2$.
Hence $J_1 = J_2$. This proves the theorem.
\end{proof}
\section{Remarks}
In case $G$ is of type $I(p)$ group and $A$ is a maximal sum-free subset of
$G$ such that $\alpha(A) > \frac{1}{3} + \frac{1}{3(p + 1)}$ then  $\alpha(A)
= \mu(G)$. For the proof of this one may see~\cite{GR}.  But in case $G$ is of
type $III$ then there exist $A$ such that $A$ is a maximal sum-free set of
cardinality $\mu(G)n - 1$.  For this consider the following example.\\ \\
\noindent{\bf Example:} $ G = (\Z/7\Z)^2$ and $A = \pi_2^{-1}\{3\} \cup (0, 2)
\cup (1, 2) \cup ( \pi_2^{-1}\{4\} \setminus \{ (0, 4), (1, 4), (2, 4) \} )$,
where $\pi_2: G \to \Z/7\Z$ is a natural projection to second co-ordinate.\\
\\ Therefore Theorem~\ref{BGCH} does not give complete characterisation of
all  large maximal sum-free subsets of $G$.\\ \\
\noindent
 In general Hamidoune and
 Plagne~\cite{Plagne1} have studied $(k, l)$ free subsets of  finite 
abelian groups. For any positive integer $t$ we define $tA$ by
\[
tA = \{ x \in G: x = \sum_{i = 1}^{t}a_i, \quad a_i \in A \quad
\forall i\}.
\]
Given any two positive integers we say $A$ is $(k, l)$ free 
if $kA \cap lA = \phi$. In case 
$k -l  \equiv 0 {\rm (mod \ m)}$, where $m$ is the exponent of $G$ then it is 
easy to check that there is no set apart from empty set which is $(k, l)$ free.
Therefore one assume that $gcd(|G|, k-l) = 1$. 
We denote the density of the largest $(k, l)$ free set by $\lambda_{k, l}$. 
Hamidoune and Plagne~\cite{Plagne1} conjectured the following.
\begin{conj}
 Let $G$ be a finite abelian group and $m$ is the exponent
of $G$. Let $k, l $ be positive 
integers  such that $gcd(|G|, k-l) = 1$. Then the density of largest $(k, l)$
free set is given by the following formula.
\[
\lambda_{k, l} = max_{d | m} \frac{[ \frac{d - 2}{k + l}] + 1}{d}
\]
\end{conj}
They~\cite{Plagne1} proved the above conjecture in  case when there exist a
 divisor $d_0$  of $m$ such that $d_0$ is not congruent to $1$ modulo $k + l$.
 In this situation they also showed that given any $(k, l)$ free set of
 density $\lambda_{k, l}$,  there exist a positive integer $d$, a surjective
 homomorphism $f: G \to \Z/d\Z$, a set $B \subset \Z/d\Z$ such that $B$ is
 $(k, l)$ free and $A = f^{-1}(B)$.
That is any $(k, l)$ free set of density $\lambda_{k, l}$ is an inverse
image of a $(k,l)$ free subset of a cyclic group. One may ask the following
 question:
\begin{Q}
 Let $G$ be a finite abelian group and $m$ is the exponent of $G$.
Suppose $k, l$ are positive integers such that all the divisors of $m$ are
 congruent to $1$ modulo $k + l$. Is it true that any $(k, l)$ free set of 
density $\lambda_{k, l}$ is an inverse image of a $(k, l)$ free set of a
 cyclic group?
\end{Q} 
We have already seen that the answer of the above question is negative in case
$k = 2$ and $l = 1$. The following arguments show that
 the answer of the above question is neagtive 
for an arbitrary value of $k$ provided $l = 1$.
\\ In case all the divisors of $m$ are
 congruent to $1$ modulo $k + l$, it is easy to check that
\[
max_{d | m} \frac{[ \frac{d - 2}{k + l}] + 1}{d} = 
\frac{ [\frac{m - 2}{k + l}] + 1}{m}.
\]
Now consider the following example.\\ \\ {\bf Example:} Let $G$ be a finite
abelian group and  $m$ is the exponent of $G$. We further assume that $G$ is
not cyclic.  Let $k $ be positive  integers  such that $gcd(|G|, k- 1) = 1$
and all the divisors of $m$ are  congruent to $1$ modulo $k + 1$.  Then $G = S
\oplus \Z/m\Z$ with $S \neq \{0\}$. Let $q = [\frac{m - 2}{k + 1}]$.  Let $x
\in \Z/m\Z$ such that  $(k -1)x \equiv 1 + q {\rm (mod \ m)}$. Let $J$ be any
proper subgroup of $S$  and
\[
A = (J \oplus \{x\}) \cup (S \oplus \{x + 1, x + 2, \cdots, x + q\} )
\cup (J^c \oplus \{ x + q + 1\}).
\]  
\\
\\
\noindent
If $A$ is in the example above then it is easy to check that $A$ is $(k,1)$
free and density of $A$ is $\frac{ [\frac{m - 2}{k + 1}] + 1}{m}$.
Hamidoune and Plagne~\cite{Plagne1} also proved the
above conjecture for all cyclic groups.
Using this and 
following the arguments as in section~\ref{secstab} it follows that 
stabiliser of $A$ is $J$. 
This shows that if $A$ is the set as in above example then $A$ is
not an inverse image of 
any $(k, 1)$ free subset of a cyclic group.\\
\\
{\bf Acknowledgment:} The authors would like to thank Dr. D. Surya Ramana
for carefully reading the manuscript and various useful comments.  
A substantial part of this work was done when the second author was a research
scholar at the Institute of Mathematical Sciences Chennai, whome he thanks for
its support.

\end{document}